\documentclass[preprint,12pt]{elsarticle}

\usepackage{hyperref} 
\usepackage{amssymb}
\usepackage{amsthm}
\usepackage{lineno}
\usepackage{amsmath}
\usepackage{epstopdf}
\usepackage{graphicx}
\usepackage{color}

\journal{Applied Mathematical Modelling}

\theoremstyle{definition}

\newtheorem{definition}{Definition}[section]
\newtheorem{remark}{Remark}[section]
\newtheorem{lemma}{Lemma}[section]

\newcommand{\REV}[1]{\textcolor{black}{#1}}

\renewcommand{\u}{\varphi}
\renewcommand{\a}{\bar\alpha}
\newcommand{\R}{\mathbb R}
\newcommand{\tr}{\textup{trace}}
\newcommand{\zmax}{z_{\textup{max}}}
\newcommand{\zmin}{z_{\textup{min}}}
\newcommand{\Gv}{\mathbf v}
\newcommand{\Gx}{\mathbf x}
\newcommand{\Gp}{\mathbf p}
\newcommand{\GX}{\mathbf X}
\newcommand{\GY}{\mathbf Y}
\newcommand{\GH}{\mathbf H}
\newcommand{\Grad}{\mathbf \nabla}
\newcommand{\Gg}{\mathbf {\hat g}}
\newcommand{\Gn}{\mathbf {\hat n}}
\newcommand{\be}{\begin{equation}}
\newcommand{\ee}{\end{equation}}

\begin{document}

%%%%%%%%%%%%%%%%%%%%%%%%%%%%%%%%%%%%%%%%%
%%%%%%%%%%%%%%%%%%%%%%%%%%%%%%%%%%%%%%%%%
\begin{frontmatter}


\title{A Level Set Based Method for Fixing \\ Overhangs in 3D Printing}

\author[label1]{Simone Cacace}
\author[label2]{Emiliano Cristiani}
\author[label3]{Leonardo Rocchi}
\address[label1]{Dipartimento di Matematica, Sapienza -- Universit\`a di Roma, Rome, Italy}
\address[label2]{Istituto per le Applicazioni del Calcolo, Consiglio Nazionale delle Ricerche, Rome, Italy (corresponding author) \texttt{e.cristiani@iac.cnr.it}}\address[label3]{School of Mathematics, University of Birmingham, Birmingham, UK}

\begin{abstract}
3D printers based on the additive manufacturing technology create objects layer-by-layer dropping fused material. As a consequence, strong overhangs cannot be printed because the new-come material does not find a suitable support over the last deposed layer. In these cases, one can add support structures (scaffolds) which make the object printable, to be removed at the end. In this paper we propose a level set based method to create object-dependent support structures, specifically conceived to reduce both the amount of additional material and the printing time. 
We also review some open problems about 3D printing which can be of interests for the mathematical community. 
\end{abstract}

\begin{keyword}
level set method \sep Hamilton-Jacobi equations \sep evolving interface \sep support structure \sep scaffolding \sep CAD software \sep additive manufacturing \sep fused deposition modelling \sep digital fabrication
\MSC[2010] 65D17 \sep 35F21 %(2000 is the default)
\end{keyword}

\end{frontmatter}

%\linenumbers

%%%%%%%%%%%%%%%%%%%%%%%%%%%%%%%%%%%%%%%%%
%%%%%%%%%%%%%%%%%%%%%%%%%%%%%%%%%%%%%%%%%
\section{Introduction}\label{sec:intro}
Is a new industrial revolution coming? Many people think so: 3D printers are able to create almost any solid object one can image and replicate existing ones. Nowadays, the price of a 3D printer is small enough to allow many people to have one at home, and create their own plastic objects. Within a decade, some products may be downloaded from the Internet for printing at home, causing a revolution in the market of such a small objects.  %One can also print at home some parts of an object and buy the missing ones at shops, then assemble the final product by himself. 
Most important, the number of printable materials is growing and it is already possible printing an object mixing different materials. We leave to futurists the comments about the time when 3D printers will be able to fully replicate themselves.

While the computer science literature about 3D printing is already rich in algorithms, optimization techniques and applications, the mathematical literature is basically null. This means that advanced mathematical tools based on PDEs, optimal control theory and variational methods are, so far, little explored. In order to fill the gap and promote the solution of the engineering issues related to CAD 3D printer software, in the next section we propose a bird-eye view over typical open problems encountered by practitioners who use 3D printers based on Fused Deposition Modeling (FDM).

\textit{Main goal.} The core of the paper is devoted to the solution to a particular problem, namely fixing the overhangs. When FDM technology is employed, the solid object is created layer by layer, starting from the lowest one. As a consequence, each layer can only be deposited \textit{on top} of an existing surface, otherwise the print material falls and solidifies ``in the air''. In doing this, little exceptions can be handled, i.e.\ the upper layer can protrude over the lower layer within a certain limit. The more the material cools down rapidly and the extruder moves slowly, the more the limit can be increased. If the overhang exceeds the hardware limit, an additional support must be necessarily added, in order to make the object printable. Note that the support structures are meant to be removed at the end of the process, and thus they represent wasted material. Even more important, they represent an additional source of printing time.

\textit{Related work.} The overhang problem was already investigated in the computer science and engineering literature, and some solutions were proposed \cite{barnett2015AM, dumas2014TOG, huang2009IJAMT, qiu2015proc, strano2013IJAMT, vanek2014CGF}. In most cases, support structures fill either densely or sparsely the free space encountered when a part is projected downward in its build orientation, see left object in Fig.\ \ref{theta_overhang}(b).
The difference between the methods is in how much material is used, the reliability of the supports, and the type of material which can be used. 
The paper \cite{barnett2015AM} proposes two support geometry algorithms particularly suitable for weak support materials. The paper \cite{dumas2014TOG} proposes an algorithm for the automatic generation of horizontal bridges and vertical pillars, connected in such a way to create a hierarchical structure. 
The paper \cite{qiu2015proc} uses a cone-based scan to detect the closest points which can serve as a support base (upon the model itself or the build plate) for any overhanging point. 
The paper \cite{huang2009IJAMT} uses instead slant hourglass-like pillars.  
The paper \cite{strano2013IJAMT} proposes to create cellular supports, riddling dense structures with holes. 
The paper \cite{vanek2014CGF} proposes an algorithm which creates thin tree-like hierarchical support structures, similar (but more efficient) to the ones generated by the software Autodesk\textregistered \ Meshmixer\textregistered \ v2.9\footnote{http://www.meshmixer.com/}.  
%ci sarebbe pure yuxin2015 ma usa un metodo parecchio diverso.

In this paper we propose to ``enlarge'' the object in such a way that supports are no longer needed. In particular, we avoid the creation of pillars which touch the build plate by means of optimally shaped chamfers, suitably placed below hanging parts, see right object in Fig.\ \ref{theta_overhang}(b).

%%%%%%%%%%%%%%%%%%%%%%%%%%%%%%%%%%%%%%%%%
%%%%%%%%%%%%%%%%%%%%%%%%%%%%%%%%%%%%%%%%%
\section{What every mathematician should know about 3D printing}\label{sec:openpbs}
In the context of 3D printers there exist several open problems and modelling needs. Generally speaking, the main issues come from the fact that software solutions are not object-dependent, and not change during the printing time, whereas each object (and each layer!) has its own peculiarities. An exhaustive bibliography is out of the scope of the paper, therefore for each problem we point out just a few significant references.

\textbf{Infill.} Printing fully solid objects is often not convenient because of the large quantity of material to be used. Shape optimization tools can give the optimal way to hollow out the object, reducing the overall material volume and keeping at the same time the desired rigidity and printable features. The problem reduces to finding the optimal inner structure supporting the whole object from the inside \cite{wang2013TOG} or partitioning the object to print hollow parts \cite{wei20163DR}. 

\textbf{Orientation \& supports.} In some cases the object is not 3D-printable due to the presence of hanging parts. In this case one should find the orientation of the solid which minimizes the hanging parts \cite{ezair2015CeG} and then the minimal amount of additional material needed to support the hanging parts. The latter problem is the one we consider in this paper. 

\textbf{Balancing.} It is important to ensure that during the printing process (and once it is finished), the object can lie in equilibrium without falling down. This problem can be solved by trying to balancing in an appropriate way the mass of the object and by creating cavities in the inner structure so that it stands in its natural pose without requiring any glue or pedestal \cite{christiansen2015CAD, prevost2013ACM}.

\textbf{Partitioning.} Sometimes it is necessary to divide a 3D model into multiple printable pieces, so as to save the space, to reduce the printing time, or to make a large model printable by small printers \cite{attene2015CGF}. This problem was attacked by means of a level set based approach similar to the one proposed here in \cite{yao2015TOG}.

\textbf{Slicing \& toolpath generation.} Creating layers from a 3D model is a crucial step in 3D printing. Usually one computes the intersection curves between the model represented by polygonal meshes and a sequence of parallel planes. However, this procedure is not trivial in case of very complicated (self-intersecting,  overlapped) objects. Moreover, once the layers are created, the exact trajectory of the nozzle must be defined. The infill pattern must be travelled in the shortest way, continuously, without halting the manufacturing process, and minimizing the jump from the end of one sub-path to the starting point of another sub-path. Interestingly, this problem can be seen as a generalized travelling salesman problem \cite{castelino2003JMS, dhanik2010IJAMT, hildebrand2013CeG, huang2013JCISE, jin2014AM}.

%\textbf{3D printing without a 3D model.} 3D-printed objects are usually originated either by CAD softwares or by 3D scans of real objects directly. In some cases, one wants to replicate a real object which cannot be 3D-scanned, so that other techniques must be employed. The typical scenario is that of the (photometric, perspective) \textit{Shape-from-Shading} problem, where one or more photographs of the same object are used to build the 3D surface corresponding to the object. The photographs are usually taken under different points of view or different light conditions.

\textbf{Shape or shading?} 3D-printed objects replicating real objects are usually made of a different (and cheaper) material with respect to the original ones. As a consequence, it is expected that the replicated object reflects light in a different manner (different albedo, different degree of Lambertianity), thus resulting in an unsatisfactory product. In some cases it can be better creating an object \emph{with different shape but which appears as the original one}. In other words, one aims at replicating the reflectance properties of an object, not its original shape \cite{lan2013TOG}.

\textbf{Oozing.} It can happen that the machine deposits too much material in some parts of the object, or the material oozes, especially when the nozzle changes direction or stays on the same point for a long time. This issue is mainly related to the temperature of the nozzle's hot end and the pressure drop because of the filament. %, see, e.g., \cite{kubicek}. 
The nozzle's temperature, the retraction of the filament and the speed of the extruder should be related to each other and optimized with respect to the printing and travelling time (i.e.\ extruder movements with and without emission of material, respectively).

\textbf{Multi-material printing.} Let us also mention the possibility to print objects with different materials simultaneously, alternating them while printing. Materials can have different reflectance properties and transparency, and, consequently, endless combinations are possible, as well as related optimization processes. Similarly, one can coat the surface with paint, thus altering the reflectance properties \cite{vidimce2013TOG}.

%%%%%%%%%%%%%%%%%%%%%%%%%%%%%%%%%%%%%%%%%
%%%%%%%%%%%%%%%%%%%%%%%%%%%%%%%%%%%%%%%%%
\section{The level set method}\label{sec:LS}

The level set method was introduced in \cite{osher1988JCP} and since then it was successfully applied in many contexts, see e.g., \cite{osherbook, sethianbook}. It allows to track Eulerianly the evolution of a $(d-1)$-dimensional surface embedded in $\R^d$ transported by a given velocity vector field $\Gv:\R^d\to\R^d$. Let us briefly recall the method in the case of $d=3$ which is of interest for our problem.

\subsection{The level set function and the Hamilton--Jacobi equation}
It is given a bounded closed surface $\Sigma_0:U\subset\R^2\to\R^3$ at initial time $t=0$. We denote by $\Sigma_t$ its (unknown) evolution under the action of $\Gv$ at time $t$ and by $\Omega_t$ the 3D domain strictly contained in $\Sigma_t$ so that $\Sigma_t = \partial \Omega_t$, for all $t\geq 0$.
The main idea of the level set method stems on the definition of a level set function %$\u(t,x,y,z):\R^+\times\R^3\to\R^4$ such that
$\u(t,x,y,z):\R^+\times\R^3\to\R$ such that
\begin{equation}
\Sigma_t=\{(x,y,z)\ :\ \u(t,x,y,z)=0\},\quad \forall \, t\geq 0.
\end{equation}
\noindent In this way the surface is recovered as the zero level set of $\u$ at any time. Initially, the function $\u$ is chosen in such a way that
\begin{equation}\label{proprietaLS}
\u(0,x,y,z)
\left\{
\begin{array}{ll}
>0, & \textrm{if } (x,y,z)\notin \overline{\Omega_0}, \\
=0, & \textrm{if } (x,y,z)\in \Sigma_0, \\
<0, & \textrm{if } (x,y,z)\in \Omega_0.
\end{array}
\right.
\end{equation}
A typical choice for $\u(0,x,y,z)$ is the signed distance function from $\Sigma_0$, although this choice does not lead to a smooth function. It is easy to prove \cite{sethianbook} that the level set function $\u$ at any later time satisfies the following Hamilton--Jacobi equation
\begin{equation}\label{advec_eq}
\partial_t \u+\Gv\cdot \Grad\u=0,\qquad t\in\R^+,\ (x,y,z)\in\R^3,
\end{equation}
with a suitable initial condition $\u(0,x,y,z)=\u_0(x,y,z)$ satisfying \eqref{proprietaLS}. Here $\Grad=(\partial_x,\partial_y,\partial_z)$ denotes the gradient with respect to the space variables. 
One of the most appealing features of the level set method is that several geometrical properties of the evolving surface can be described by means of its level set function $\u$. For example, it possible to write the unit exterior normal $\Gn$ and the (mean) curvature $\kappa$ in terms of $\u$ and its derivatives. More precisely, we have
$$\Gn=\frac{\Grad\u}{|\Grad\u|} \quad \text{ and } \quad \kappa= \Grad \cdot \Gn.$$
If the vector field has the form $\Gv = v\Gn$ for some scalar function $v$, the equation \eqref{advec_eq} turns into
\begin{equation}
\label{LSeq}
\partial_t \u+v|\Grad\u|=0,\qquad t\in\R^+,\ (x,y,z)\in\R^3.
\end{equation}
%\noindent Both \eqref{advec_eq} and \eqref{LSeq} are first order evolutive Hamilton-Jacobi type equations of the form $\partial_t \u+H(x,y,z,D\u)=0$, where $H$ represents the Hamiltonian, say $H = \Gv \cdot D\u$ and $H = v|D\u|$ respectively. Classical analysis of these equations led to the appropriate definition of viscosity solution \cite{crandallLions1983TAMS} for which numerical consistent monotone schemes have been developed \cite{crandallLions1984TAMS}.

\subsection{Computation of the signed distance function}\label{sec:LS.df}
The computation of the signed distance function $\u_0$ is a problem \emph{per se}. In our case, we can assume that the surface $\Sigma_0$ of the object to be printed is watertight and that it is given by means of a triangulation (typically in the form of a .STL file). 
Each triangle (facet) $f$ is characterized by the 3D coordinates of its three vertices. Moreover, vertices are oriented in order to distinguish the internal and the external side of the facet.

Given a point $(x,y,z)\in\R^3$, it is easy to find the distance $d((x,y,z),f)$ between the point and the facet, so that the unsigned distance from the surface is given by
$$
d((x,y,z),\Sigma_0)=\min_f d((x,y,z),f).
$$ 
The computation of the distance's sign is more tricky since one has to check if the point is internal or external to the surface. Several methods can be employed here. For example, one can note that the solid angle subtended by the whole surface at a given point is maximal and equal to $4\pi$ iff the point is internal. Then, one can sum all the solid angles subtended by the facets at the point and check if it equals $4\pi$. If this is the case, the point is internal to the surface, otherwise it is external.

Note that the solid angle itself should be \emph{signed}, in the sense that it must be positive if the point looks at the internal part of the facet, negative otherwise. A nice algorithm to compute the signed solid angle between a point and a triangle was given by van Oosterom and Strackee \cite{oosterom1983}.

%%%%%%%%%%%%%%%%%%%%%%%%%%%%%%%%%%%%%%%%%
%%%%%%%%%%%%%%%%%%%%%%%%%%%%%%%%%%%%%%%%%
\section{Fixing overhangs}\label{sec:fixingoverhangs}
In this section we propose a solution for fixing the overhang issue in 3D printing and, in most cases, getting rid of long support structures extended until the build plate. We want to use the level set method by considering the domain $\Omega$ as the object to be printed and its surface $\Sigma$ as an evolving front. Therefore, the idea is to modify the initial unprintable object $\Omega_0$ letting it evolve by an \emph{ad hoc} vector field $\Gv$ until it becomes fully printable, meaning that there are no more unprintable hanging parts. The final object $\Omega_*$ is then actually printed and the difference $\Omega_*\backslash\Omega_0$ is finally removed. Note that the difference $\Omega_*\backslash\Omega_0$ can be easily identified by standard techniques and consequently printed with a different material (e.g., a soluble filament) or with a different printing resolution.

It is useful to divide the surface $\Sigma$ of the object $\Omega$ in three subsets, on the basis of their \emph{printability}. To this end, we denote by $\Gg=(0,0,-1)$ the unit gravity vector, and again by $\Gn(x,y,z)$ the exterior unit normal to the surface $\Sigma$ of the object $\Omega$ at the point $(x,y,z)$. Moreover, let
\be
\theta(\Gn) := \arccos \big(\Gg\cdot\Gn\big) %\in [\a,2\pi-\a],
\ee
be the angle between $\Gg$ and $\Gn$.

\begin{figure}[h!]
\centering
\begin{tabular}{cc}
\includegraphics[scale = 0.33]{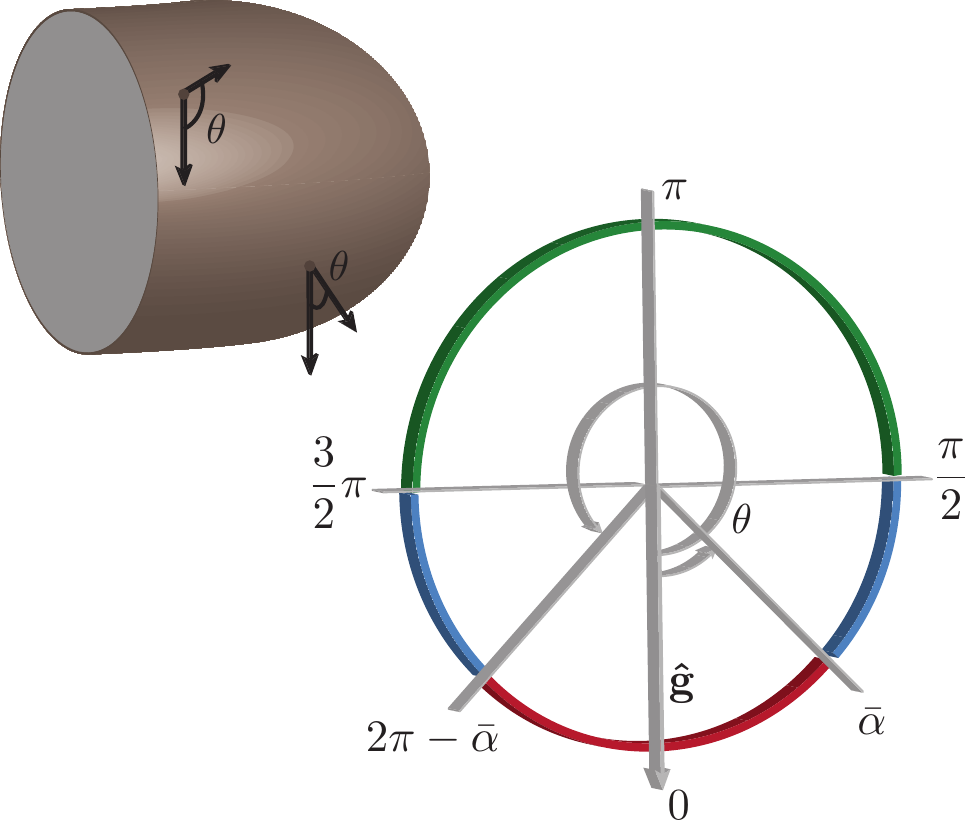} \qquad \qquad &
\includegraphics[scale = 0.24]{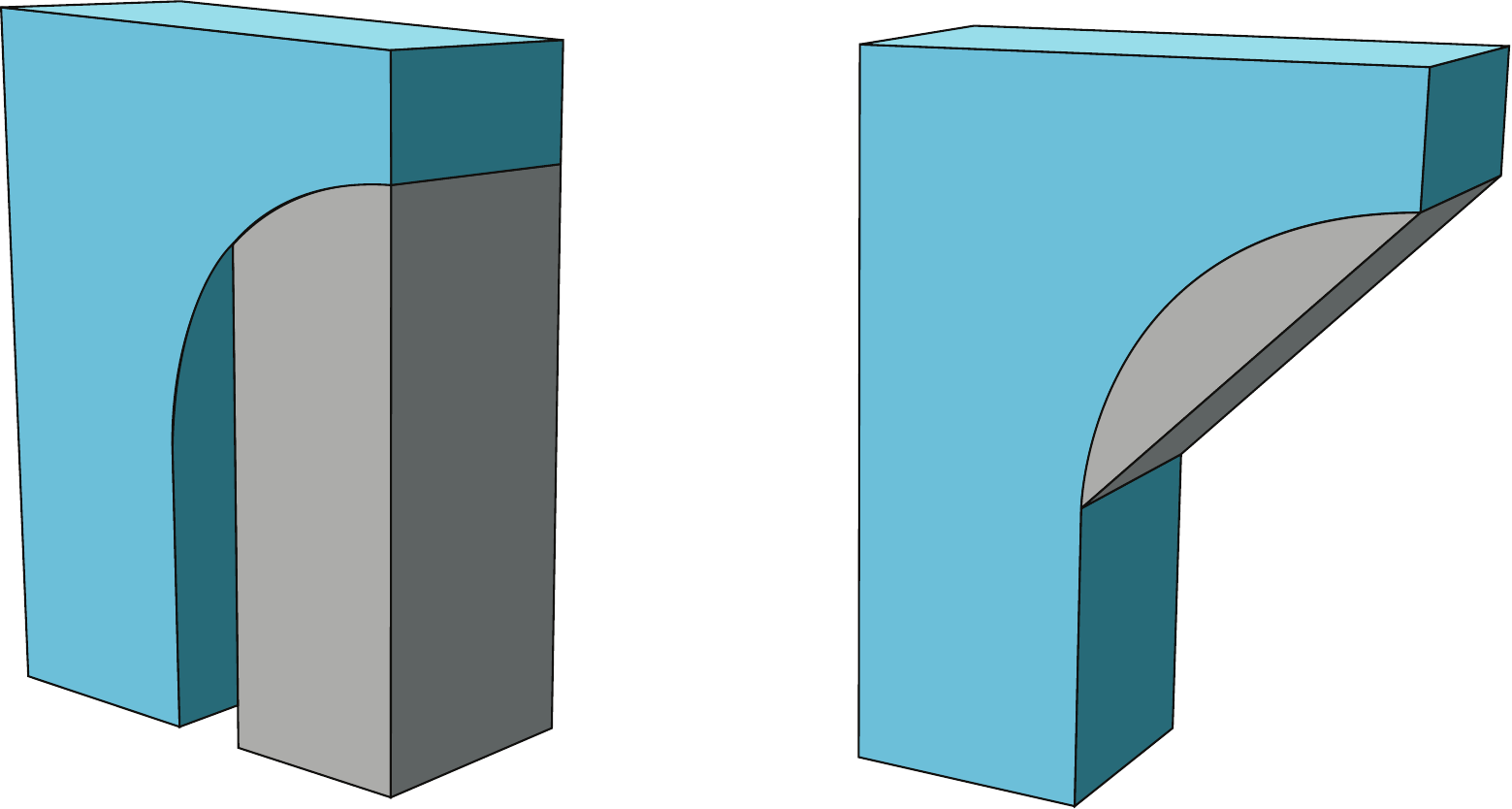} \\
\footnotesize (a) & \footnotesize (b) 
\end{tabular}
\caption{\footnotesize{(a) Unprintable (red), modifiable (blue) and safe (green) points with respect to the counter-clockwise angle $\theta$ between the gravity $\Gg$ and the normal $\Gn$. Modifiable and safe points are printable. (b) The left grey support wastes a lot of material contrary to the chamfer on the right that saves more material and keeps the printability of the overhang as well.}}
\label{theta_overhang}
\end{figure}

\begin{definition}
\label{def:printable}
\emph{A point $(x,y,z)$ of the surface $\Sigma$ is said to be}
\begin{center}
\begin{tabular}{rl}
\textit{unprintable,}  & if $\theta\in[0,\a)\cup(2\pi-\a,2\pi]$, \\
\textit{safe,}         & if $\theta\in[\pi/2,3\pi/2]$, \\
\textit{modifiable,}   & otherwise,
\end{tabular}
\end{center}
\noindent \emph{where $\a$ is a given limit angle\footnote{Typically $\a=\frac{\pi}{4}$, because of the so-called $45$ degree rule, though it actually depends on the 3D printer settings, print material, cooling, etc.}}, see Fig.\ \ref{theta_overhang}(a).
\end{definition}  

While the first two definitions are immediately clear, it is worth to spend some words on the third one. Modifiable points are indeed printable since the overhang is sufficiently small. On the other hand, it could be convenient to move those points as well in order to make printable the unprintable ones. This guarantees a sufficient flexibility to shape the object conveniently and not to create long supports like the one depicted on the left of Fig.\ \ref{theta_overhang}(b). We can extend Definition \ref{def:printable} by saying that the set of both modifiable and safe point constitute the overall printable points. %, see Fig.\ \ref{theta_overhang}(a).

The rest of the section will be devoted to the construction of the vector field $\Gv$. A suitable choice is the one used in equation \eqref{LSeq} where $\Gv=v\Gn$ for some scalar function $v$, possibly depending on $\Gn$ and $\kappa$. 

In the following we denote by 
$$P(\omega):=\omega^+ \quad \text{and} \quad  M(\omega):=\omega^-, \qquad \omega\in\R,$$ 
the positive and negative part, respectively.

\medskip

\noindent \textbf{Positivity and build plate.} We need to be guaranteed that $\Omega_0\subseteq \Omega_t\subseteq\Omega_*$, for all $t\geq 0$, since once the object is printed we can remove material but not add new one. This is why we need $v\geq 0$, i.e.\ the movement of each point of the surface $\Sigma$ has to be along the normal exterior direction $\Gn$. Furthermore, the object cannot move under the build plate, supposed at a fixed $z=\zmin \in \R$. Then we impose $v=0$ if $z \leq \zmin$.

\noindent \textbf{Movement of unprintable points.} We introduce the term
\be
v_1(\Gn;\a):=P(\cos\theta(\Gn)-\cos\a),
\ee
which lets the unprintable points move outward. The speed is higher whenever $\theta$ is close to $0$, which represents the (hardest) case of a horizontal hanging part.

\noindent \textbf{Rotation.} It is convenient introducing a rotational effect in the evolution which avoids the unprintable regions to evolve ``as it is'' until they touch the build plate. To this end we introduce the term $(\zmax-z)$, where $\zmax \in \R$ is the maximal height reached by the object. This term simply increases the speed of lower points with respect to higher ones. This makes the lower parts be resolved (or eventually touch the built plate) before the higher parts, thus saving material. 

\noindent \textbf{Movement of modifiable points.} Modifiable points are moved, if necessary, by means of the following term in the vector field
\be
v_2(\kappa) := M(\kappa).
\ee
It moves outward the points with negative curvature until it vanishes, i.e.\ the surface is locally flat. In particular, it moves concave corners and let modifiable points become a suitable support for the still unprintable points above.

\noindent \textbf{Blockage of safe points.} Finally, it is necessary to exclude from the evolution the safe points of the object. In order to identify them, we use the sign of the third component $n_3$ of the unit exterior normal vector $\Gn$. %More precisely, we use the term $n_3^-(x,y,z)$ in order to block safe points.

\medskip

By putting together all the terms we end up with 
%$\Gv = v(x,y,z,\Gn,\kappa;\a)\Gn$, where the velocity normal function reads as 
%\begin{equation}
%\label{complete_v}
%v(x,y,z,\Gn,\kappa;\a) := 
%\left\{
%\begin{array}{ll}
%\chi(z)\Big( A\, v_1(x,y,z,\a) (\zmax-z) + B\,v_2(x,y,z)\Big) & \text{ if } n_3<0 \\
%0 & \text{ if } n_3\geq 0
%\end{array}
%\right.
%\end{equation}
\begin{multline}\label{complete_v}
v(x,y,z,\Gn,\kappa;\a) :=\\ 
\left\{
\begin{array}{ll}
C_1\, (\zmax-z)v_1(\Gn;\a) + C_2\,v_2(\kappa), & \text{ if } n_3<0 \text{ and } z > \zmin,\\
0, & \text{ otherwise, }
\end{array}
\right.
\end{multline}
with $C_1, C_2>0$ positive constants (model parameters). %and $\chi(z)=1$ if $z > \zmin$ and $\chi(z)=0$ otherwise. 
%\[
%\chi(z) := 
%\begin{cases}
%1 & \text{if $z > \zmin$}     \\
%0 & \text{if $z \leq \zmin$}.
%\end{cases}
%\]
The result expected from a such vector field is an evolution similar to the one depicted on the right in Fig.\ \ref{theta_overhang}(b), corresponding to a support whereby the angle $\theta$ in each of its point is less or equal to $\a$.

\medskip

The surface evolution relative to equation \eqref{LSeq} must be stopped at some final time $T>0$. Rather than waiting that the velocity field vanishes completely, it is convenient to check directly (at every time $t<T$) whether the overall surface is printable or not, according to Definition \ref{def:printable}. More precisely, we stop the evolution when all the points belonging to the zero level set are safe or modifiable, i.e., printable.

\begin{remark}\label{rem:optimality}
(Optimality of the final surface)
By construction, the surface always evolves towards a printable object. 
Indeed, any non-printable part of the surface is forced to move downward, and the surface has to stop once the build plate is reached.
Nevertheless, we have no guarantee that the final object is ``optimal'' in terms of additional printing material. In the worst-case scenario the surface evolves until it touches the build plate, obtaining something similar to the results depicted on the left in Fig.\ \ref{theta_overhang}(b). For instance, this is the case of a perfectly symmetric bridge-shaped object, unless some symmetry-breaking terms are added in the evolution model. Likely, the method works fine in most cases, as one can see in section \ref{sec:numericaltests}, where several objects are tested. 
\end{remark}

%%%%%%%%%%%%%%%%%%%%%%%%%%%%%%%%%%%%%%%%%
%%%%%%%%%%%%%%%%%%%%%%%%%%%%%%%%%%%%%%%%%
\section{Theoretical analysis}\label{sec:theory}
In this section we show that a slightly regularised version of the Hamilton--Jacobi equation \eqref{LSeq} with velocity field \eqref{complete_v} fits the classical theory of viscosity solutions and it is then well-posed. 

Consider a general second order PDE of the form
\begin{equation}\label{2orderHJstandardform}
\u_t+F(t,\Gx,\u,\Grad\u,\GH\u)=0,\qquad t>0,\quad \Gx\in\R^n,
\end{equation}
where $\GH\u$ is the Hessian matrix of $\u$ and $F:\R^+\times\R^n\times\R\times\R^n\times\mathcal S_n\to\R$ is continuous and $\mathcal S_n$ is the set of symmetric $n\times n$ matrices. Resorting to classical results \cite{crandall1992BAMS}, we can say that the IVP for \eqref{2orderHJstandardform} is well-posed if the function $F$ is \emph{proper} for any fixed $t\in[0,T]$, i.e.
\begin{align}\label{def:proper}
\forall t \quad F(t,\Gx,r,\Gp,\GX)\leq F(t,\Gx,s,\Gp,\GY)\quad \text{whenever } r \leq s \text{ and }  \\ \nonumber \GX-\GY \text{ is positive semi-definite}.
\end{align}
Before writing our equation in the form \eqref{2orderHJstandardform}, let us note that $M(\omega)=P(-\omega)$ for all $\omega\in\R$, and $M(c\omega)=cM(\omega)$, for all $\omega\in\R$ and $c>0$. 
Moreover, let $H:\R\to\{0,1\}$ be the Heaviside function and state the following equality:
\begin{lemma}\label{lemma}
Given any function $u\in C^2(\R^n;\R)$, we have
$$
\textup{div}\left(\frac{\Grad u}{|\Grad u|}\right)| \Grad u|=
\tr\left(\left(\mathbf I-\frac{\Grad u \otimes \Grad u}{|\Grad u|^2}\right)\GH u\right),
$$
where $\mathbf I$ is the $n\times n$ identity matrix and $(\mathbf{a}\otimes \mathbf{b})_{i,j}=a_ib_j$ for all $\mathbf{a}, \mathbf{b}\in\R^n$ and $i,j=1,\ldots,n$.
\end{lemma}
The proof of the Lemma is postponed in the Appendix.

Making explicit the dependence on $\u$, we can rewrite the speed as
\begin{multline*}
v(z,\nabla\u,\GH\u;\a)=\\ 
\left[C_1(\zmax-z)M\left(\frac{\partial_z\u}{|\nabla\u|}+\cos\a\right)+C_2 M\left(\textup{div}\left(\frac{\nabla\u}{|\nabla\u|}\right)\right)\right]H(-\partial_z\u)H(z-\zmin).
\end{multline*}
Note that in our case the function $v$ does not depend explicitly on $t$ and $\u$, but depends implicitly on $\GH\u$ by means of the divergence operator. 
In our case we have $F=v|\nabla\u|$, and then, using Lemma \ref{lemma}, we have
\begin{multline*}
F(z,\nabla\u,\GH\u;\a)=\\
\Bigg[C_1(\zmax-z)M\left(\partial_z\u+|\nabla\u|\cos\a\right)+ \hskip5cm \\ 
\hskip4cm C_2 M\left(\textup{div}\left(\frac{\nabla\u}{|\nabla\u|}\right)|\nabla\u|\right)\Bigg]H(-\partial_z\u)H(z-\zmin)=\\
\Bigg[C_1(\zmax-z)M\left(\partial_z\u+|\nabla\u|\cos\a\right)+\hskip5cm \\ 
\hskip 2cm C_2 M\left(\tr\left(\left(\mathbf I-\frac{\nabla\u \otimes \nabla\u}{|\nabla\u|^2}\right)\GH\u\right)\right)\Bigg]H(-\partial_z\u)H(z-\zmin).
\end{multline*}
Following again \cite{crandall1992BAMS} (Example 1.2), and considering the sign of the function $M$, we are left to prove that the matrix  $\left(\mathbf I-\frac{\Gp \otimes \Gp}{|\Gp|^2}\right)$, $\Gp\in\R^3$, is positive semi-definite. A straightforward computation shows that, for all vectors $(x_1, x_2, x_3)\neq(0,0,0)$,
\begin{multline*}
\begin{array}{ccc} (x_1\!\!\! & x_2\!\!\! & x_3) \end{array}
\left(\mathbf I-\frac{\Gp \otimes \Gp}{|\Gp|^2}\right)
\left(\begin{array}{c} x_1 \\ x_2 \\ x_3 \end{array}\right)=\\ 
\frac{1}{|\Gp|}\left[(x_1p_2-x_2p_1)^2
+(x_1p_3-x_3p_1)^2
+(x_2p_3-x_3p_2)^2
\right]\geq 0.
\end{multline*}
This proves that $F$ is proper. In order to entirely fit the theoretical framework we should guarantee the continuity of the function $F$, although this is not expected to be a crucial point from the numerical point of view, since continuity cannot be actually satisfied a discrete level. An easy solution is the mollification of the Heaviside function by convolution, which makes $F$ be continuous.

%%%%%%%%%%%%%%%%%%%%%%%%%%%%%%%%%%%%%%%%%
%%%%%%%%%%%%%%%%%%%%%%%%%%%%%%%%%%%%%%%%%
\section{Numerical tests}\label{sec:numericaltests}
We solve equation \eqref{LSeq} with velocity \eqref{complete_v} by using a monotone upwind scheme based on finite differences as described in \cite[Sect.\ 6.4]{sethianbook}, with an adaptive time step in order to strictly satisfy the CFL condition. 

\medskip

\begin{figure}[h!]
\centering
\begin{tabular}{cccc}
\includegraphics[width = 0.205\textwidth]{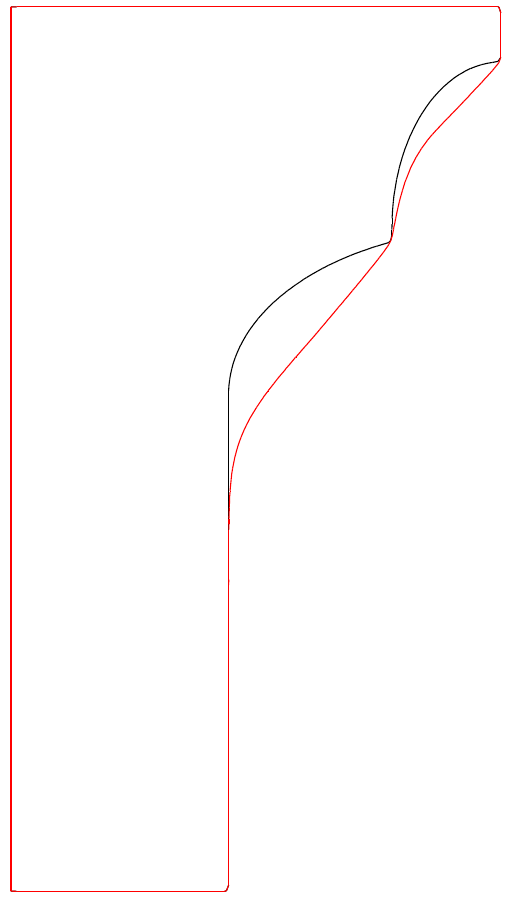} &
\includegraphics[width = 0.211\textwidth]{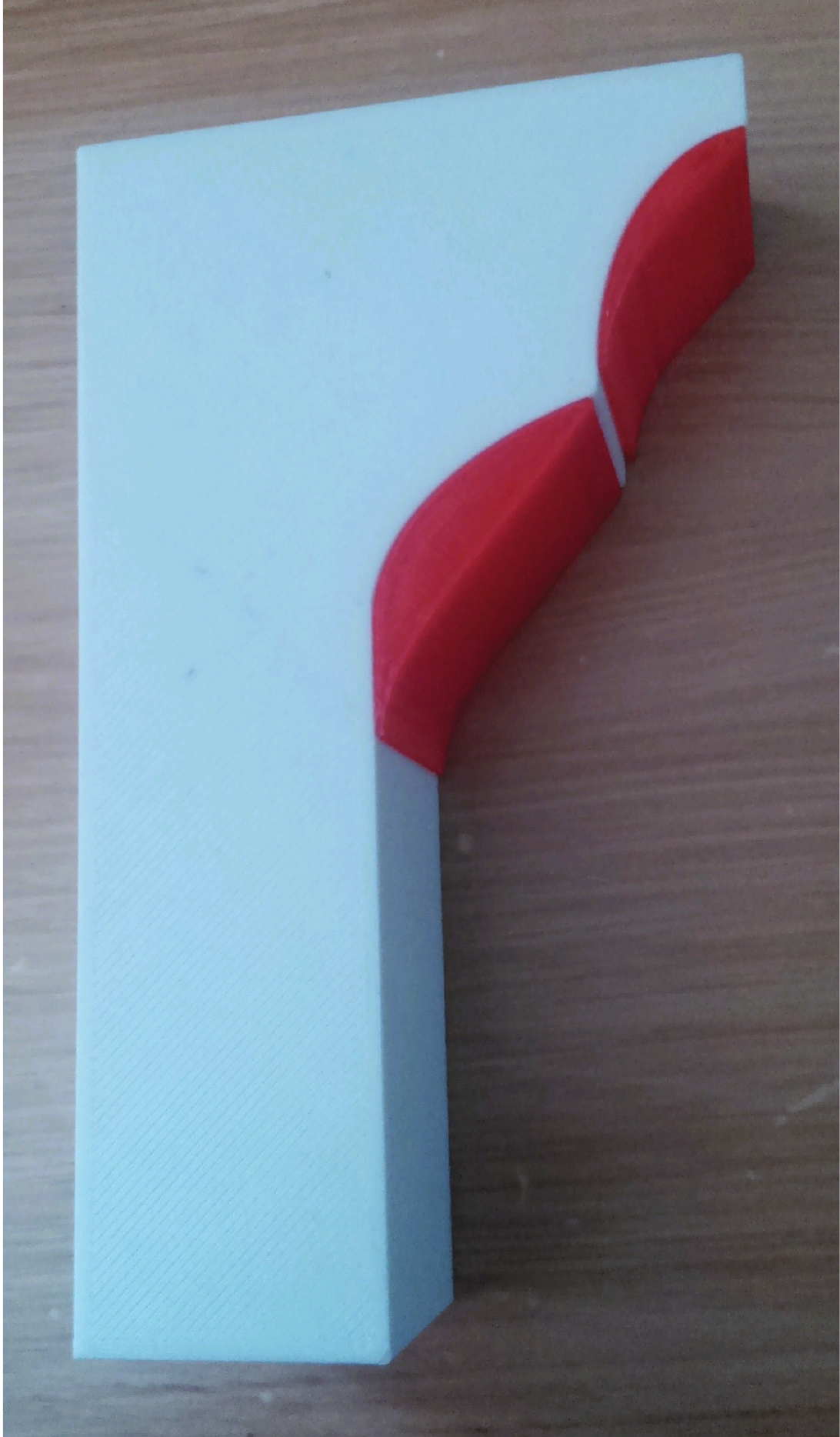} &
\includegraphics[width = 0.25\textwidth]{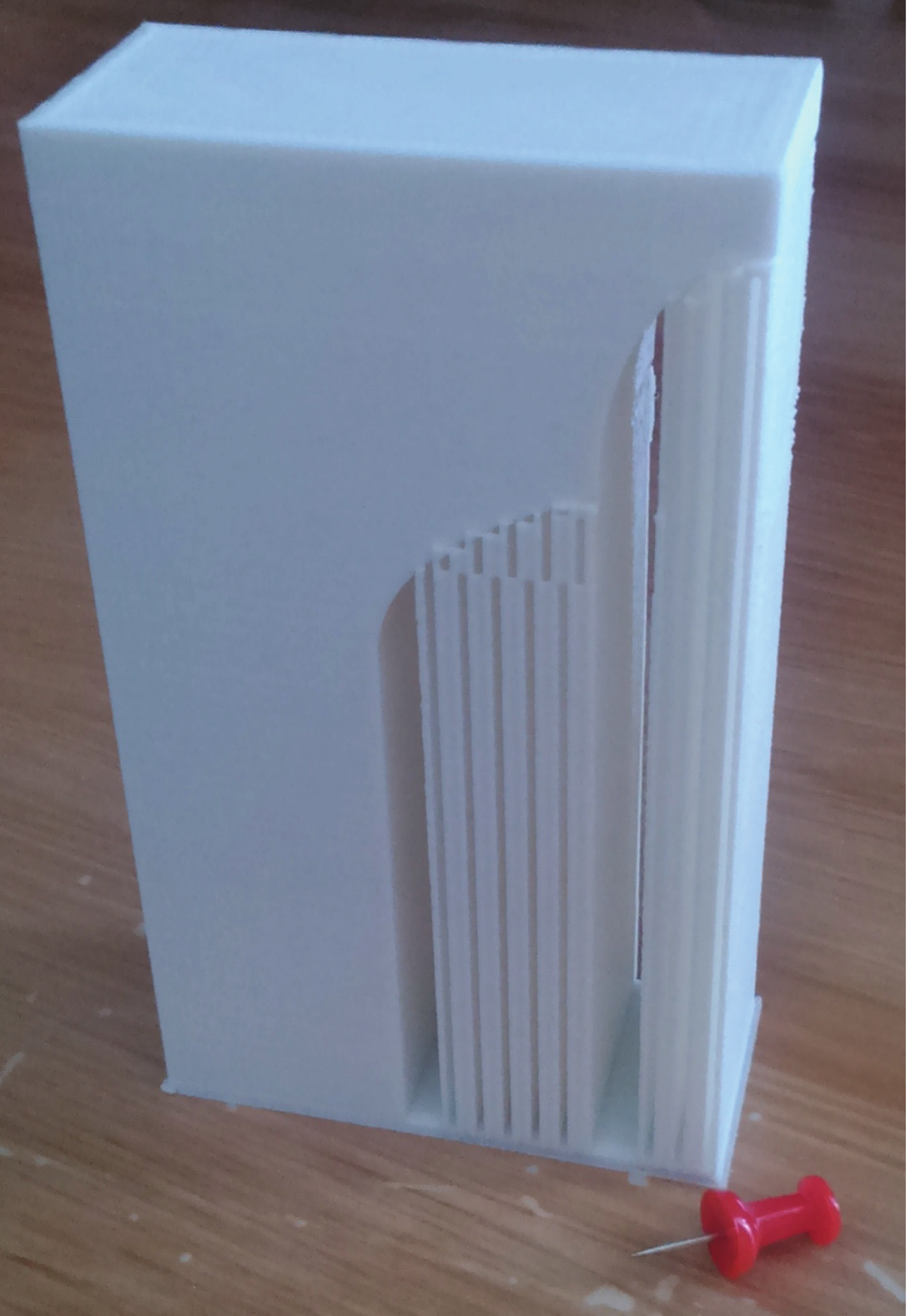} &
\includegraphics[width = 0.23\textwidth]{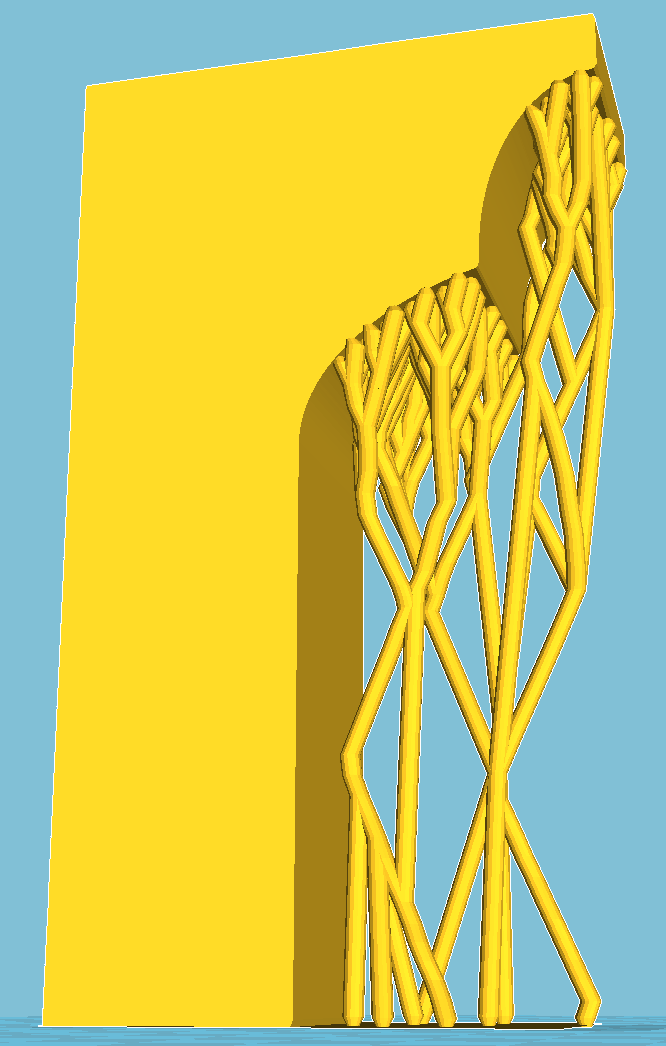} \\
\footnotesize (a) & \footnotesize (b) & \footnotesize (c) & \footnotesize (d)
\end{tabular}
\caption{\footnotesize{2D test. 
(a) Initial contour $\Sigma_0$ (black) and the optimized one $\Sigma_*$ (red) after the evolution. 
(b) Printed model with proposed support.
(c) Printed model with support structure generated by free software Cura v15.04.2. 
(d) Tree-like support structure generated by Autodesk\textregistered \ Meshmixer\textregistered \ v2.9. 
}} 
\label{2Dex}
\end{figure}

As a preliminary test, we solved a dimension-reduced problem by considering a 2D interface~$\Sigma$ with two hanging parts as the zero level set of a specific level set function $\u:\R^+\times \R^2 \to \R$. 
%We then stopped the simulation as soon as all the points of $\Gamma$ become printable according to Definition \ref{def:printable}. 
The computational domain is $[0,6]\times[0,10]$, divided in $120\times 200$ regular grid nodes. Parameters are $C_1=6$ and $C_2=0.4$. Initial and final shapes of the interface are shown in Fig.\ \ref{2Dex}(a). Moreover, by ``extruding'' the 2D domain $\Omega$, as it was a section of a real 3D object, we printed it out with the supports created from our method (Fig.\ \ref{2Dex}(b)), and keeping the scaffolding structure created by the commercial software Cura v15.04.2 (Fig.\ \ref{2Dex}(c)). 
Finally, Fig.\ \ref{2Dex}(d) shows the support structure generated by the commercial software Autodesk\textregistered \ Meshmixer\textregistered \ v2.9.

\begin{figure}[h!]
\centering
\begin{tabular}{cc}
(a) \includegraphics[scale = 0.49]{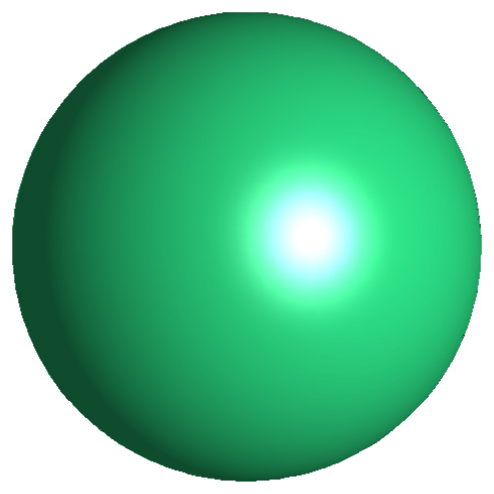} \qquad\qquad\qquad &
\includegraphics[scale=0.49]{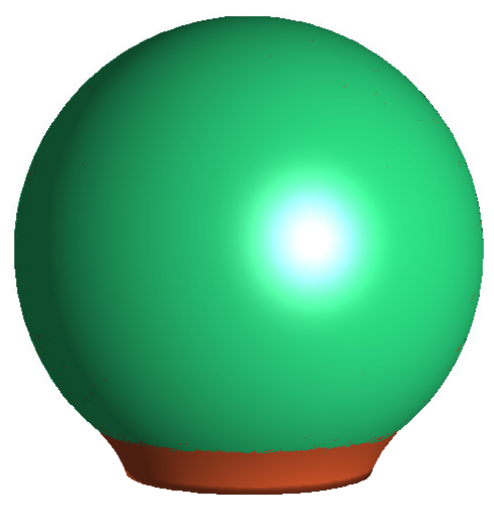} \\ \\
(b) \includegraphics[scale = 0.49]{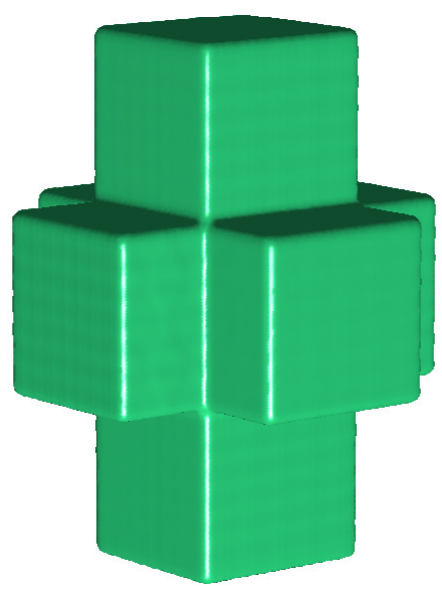} \qquad\qquad\qquad &
\includegraphics[scale=0.49]{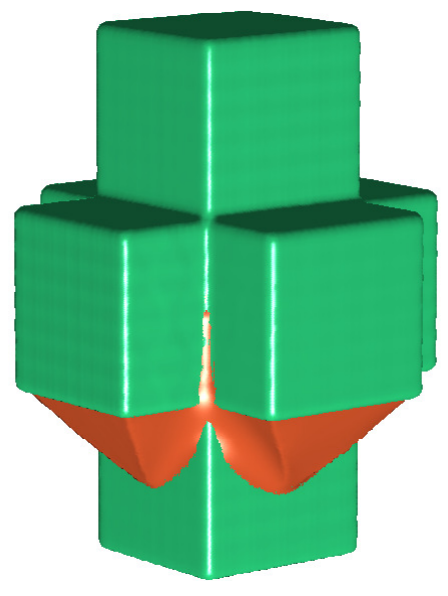}
\end{tabular}
\caption{\footnotesize{Two simple objects. Initial surface (left) and final result (right).}}
\label{fig:spherecross}
\end{figure}

Moving to real 3D problems, we tested eight objects. 
In all cases the computational domain $[-2,2]^3$ is divided in $100^3$ regular grid nodes.  
The first two examples, a sphere and a cross, are shown in Fig.\ \ref{fig:spherecross}(a,b). They have been easily obtained as the zero level set of a corresponding hyper-surface embedded in $\R^4$ and no .STL files have been required. 
The parameters used for the evolution are $C_1=0.7$, $C_2=0.3$ for the sphere and $C_1=1.5$, $C_2=0.5$ for the cross.
%%%%%%%%%%%%%%%%%%%%%%%%%%%%%%%%
%%%%%%%%%%%%%%%%%%%%%%%%%%%%%%%%

These simple numerical tests clearly show the advantage of the proposed approach: the additional material used to make the object printable is rather minimal, being concentrated in the critical zones. No evident waste of material is visible. More precisely, we see that the additional material is limited to the quantity needed to support overhangs within the maximum allowed slope.

\medskip

\begin{figure}[h!]
\centering
\begin{tabular} {c@{\hskip 0.7in}c}
\hskip-0.4cm (a) \hskip0.4cm \includegraphics[scale = 0.5]{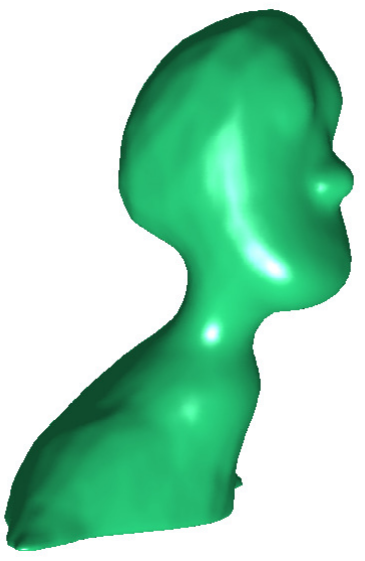}   &
\includegraphics[scale = 0.5]{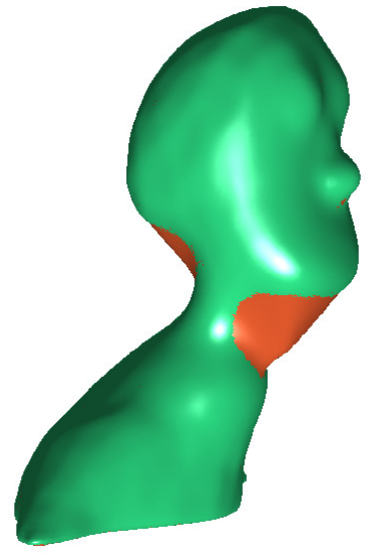} \\
(b) \includegraphics[scale = 0.5]{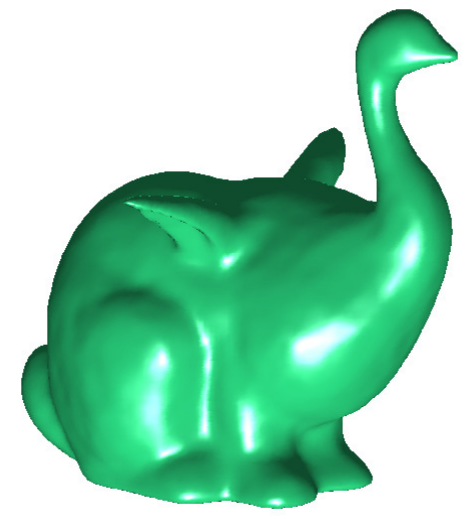}   &
\includegraphics[scale = 0.5]{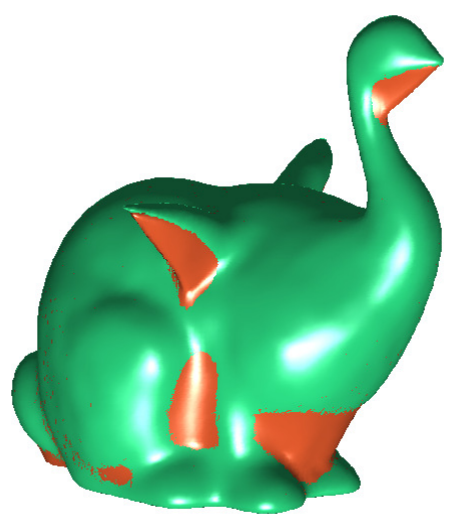} \\ \\
\ (c) \includegraphics[scale = 0.5]{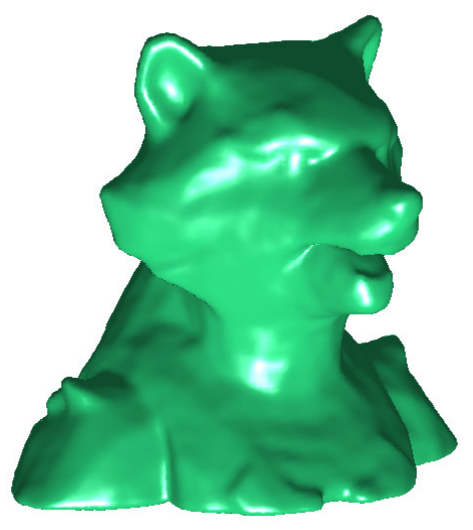} &
\includegraphics[scale = 0.5]{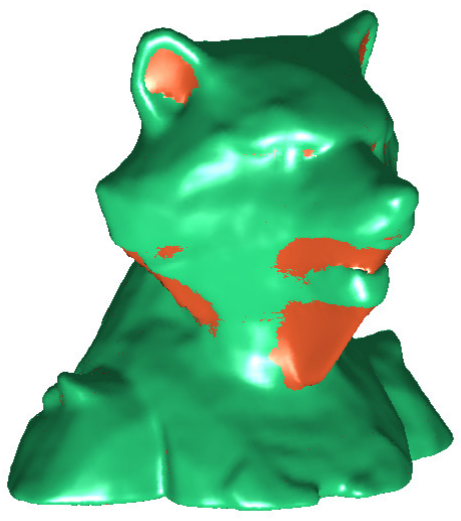} 
\end{tabular}
\caption{\footnotesize{Three objects. Initial surface (left) and final result (right).
}}
\label{fig:newtests}
\end{figure}

The following three examples are shown in Fig.\ \ref{fig:newtests}(a,b,c). 
In this case, we obtained the objects as the zero level set of the distance function (see section \ref{sec:LS.df})
for the corresponding surface given as .STL file.
The parameters used for the evolution are (a) $C_1=0.8$, $C_2=0.8$, (b) $C_1=1.8$, $C_2=0.6$, and (c) $C_1=0.6 $, $C_2=1.2$.
%%%%%%%%%%%%%%%%%%%%%%%%%%%%%%%%

Again we see that the additional material is rather minimal and concentrated in the critical zones.

\medskip

\begin{figure}[h!]
\centering
\begin{tabular} {c@{\hskip 0.7in}c}
(a) \includegraphics[scale = 0.2]{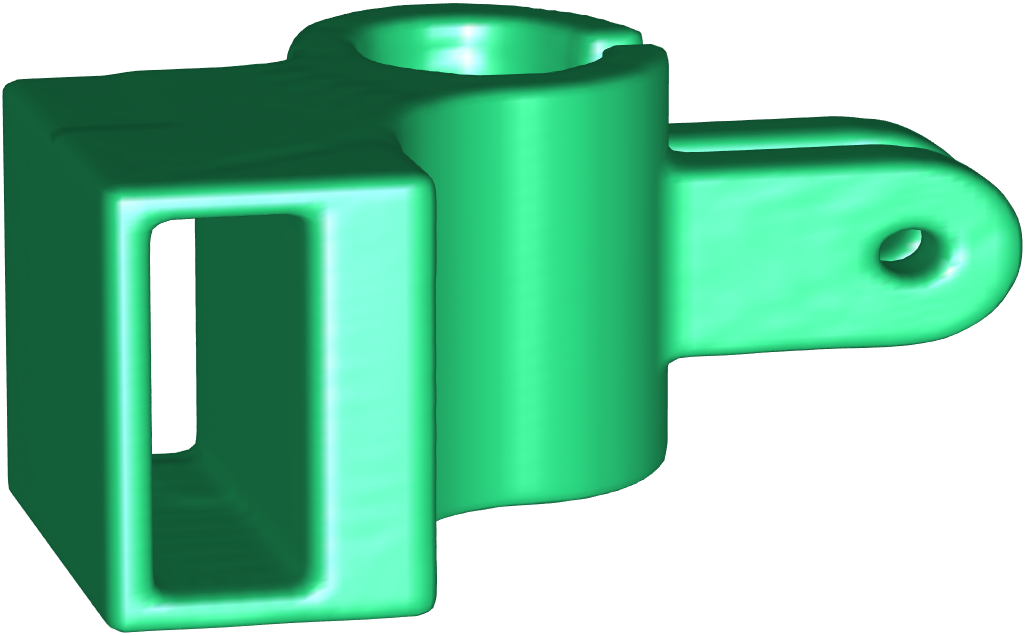}   &
\includegraphics[scale = 0.2]{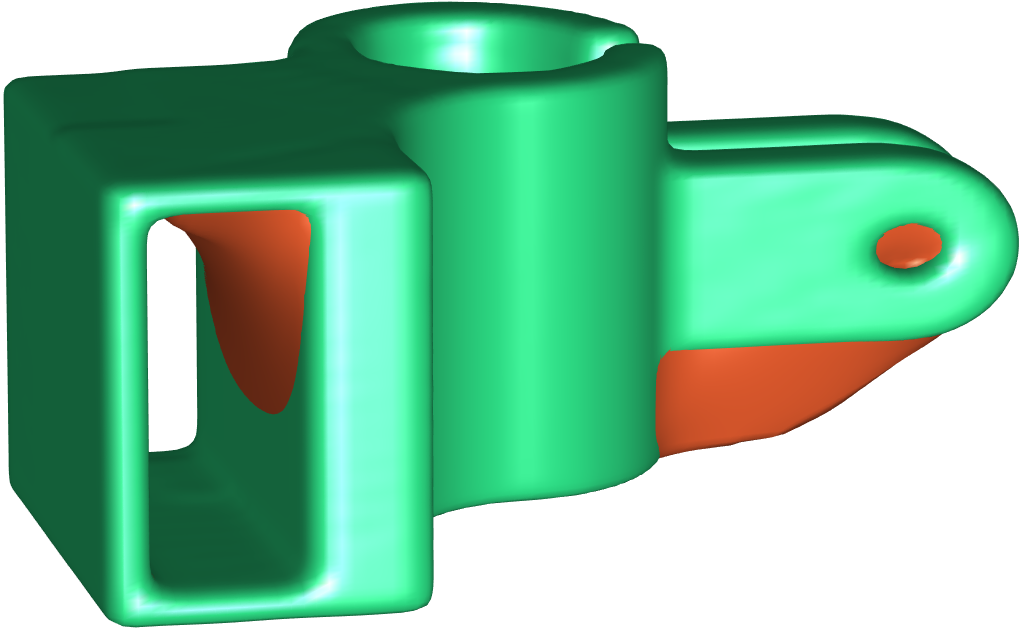} \\
(b) \includegraphics[scale = 0.2]{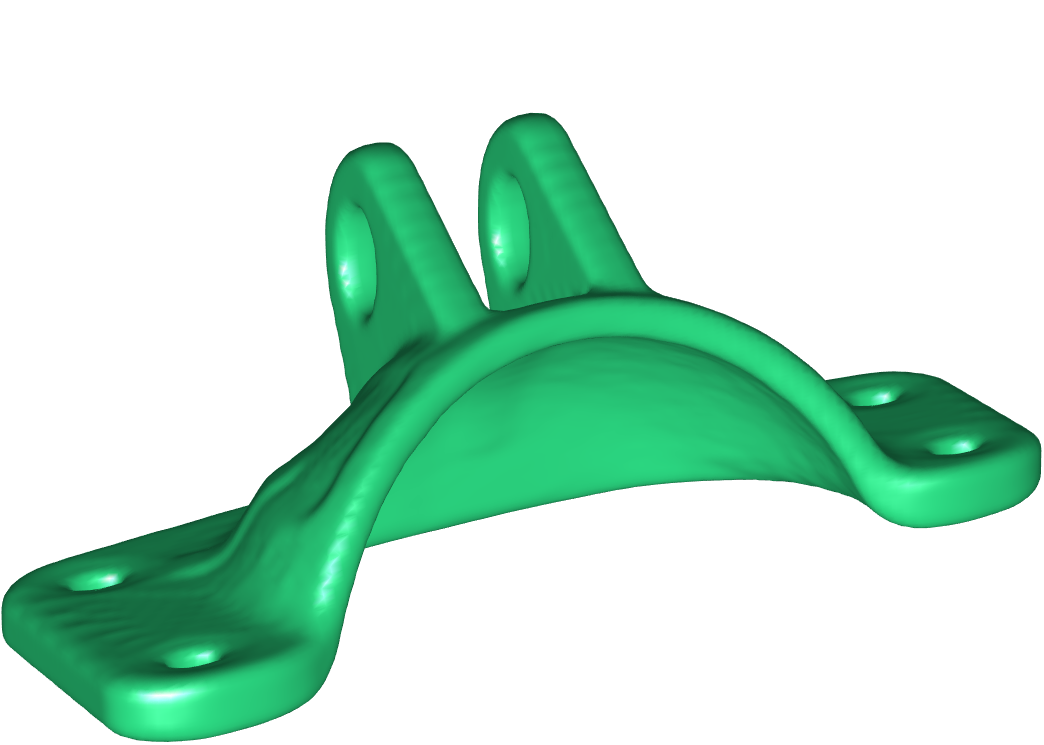}  &
\includegraphics[scale = 0.23]{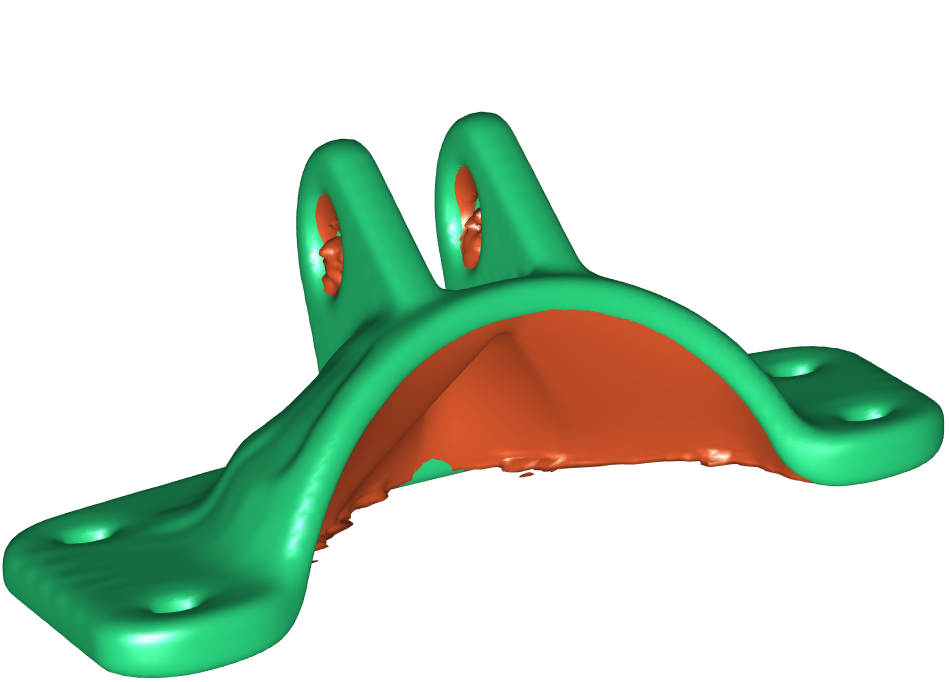}  \\ \\
\hskip-0.2cm (c) \hskip0.2cm \includegraphics[scale = 0.2]{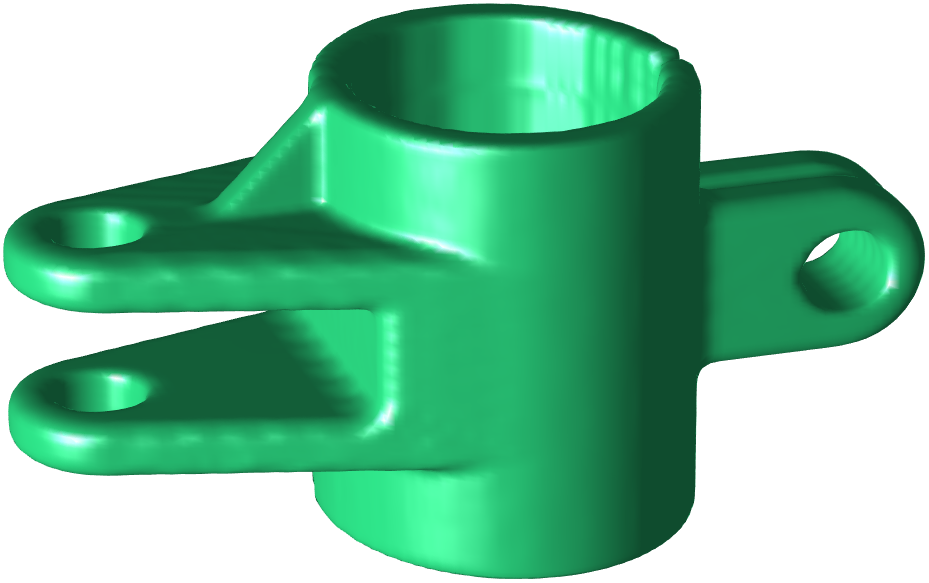} &
\includegraphics[scale = 0.15]{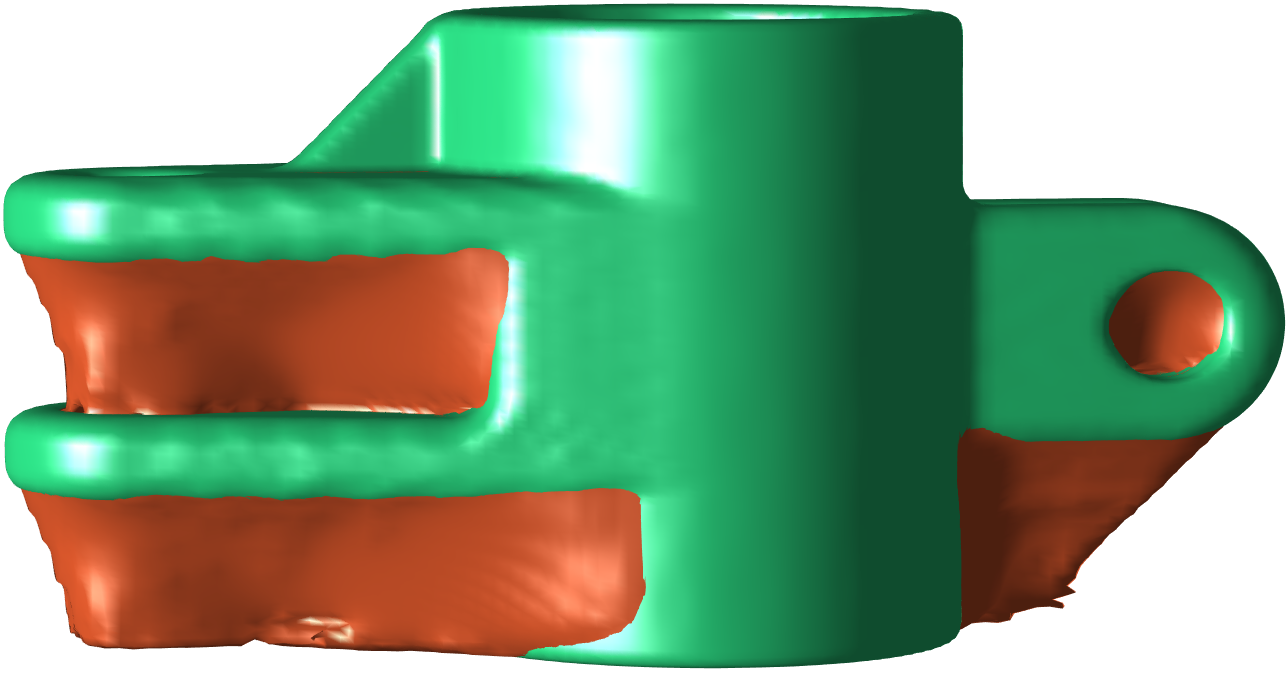} 
\end{tabular}
\caption{\footnotesize{Three brackets. Initial surface (left) and final result (right).
}}
\label{fig:brackets}
\end{figure}

\REV{The last three examples are shown in Fig.\ \ref{fig:brackets}(a,b,c). 
In this case we tried to fix overhangs of some mechanical components, starting again from the corresponding .STL files.
The parameters used for the evolution are 
(a) $C_1=3.5$, $C_2=0.9$, 
(b) $C_1=1.2$, $C_2=1.6$, and 
(c) $C_1=5.6$, $C_2=0.3$.}
%%%%%%%%%%%%%%%%%%%%%%%%%%%%%%%%

\REV{In the first two cases the evolution is close to the optimal one since no evident waste of material is visible. In the last case instead, the horizontal holes on the overhanging part are a challenge for the proposed method. They hinder the rotation and make the surface evolve vertically until the build plate is reached, thus realizing the worst-case scenario discussed in Remark \ref{rem:optimality}.}

\section*{Conclusions and future work}
We have introduced a level set based method to create \emph{ad hoc} chamfers in additive manufacturing, avoiding in most cases the use of classical vertical support structures. Moreover, in the worst-case scenario the evolved surface will not be worse than the one obtained with commercial software.

The main drawback of the proposed approach is that objects with small or sharp details clearly require a quite fine computational grid, thus rising the CPU time. Note that this is not a limitation of the evolution model, rather a limitation of the level set method itself.

Let us also stress here that we consider only the overhanging issue in printability, although the use of the support structure is not just for overhangs. As recalled in section \ref{sec:openpbs}, it also keeps the whole model from falling because of the gravity. Moreover, for printers that use the same material for both the model and the support structure, removing the support structure from the model can become difficult too. For all these reasons the proposed method is more suitable for printers that use different materials for the support structure, such as polyjet printers.

We hope that this study can pave the way to shape optimization methods based on the coupling of the level set method and the shape derivatives \cite{burger2005EJAM}. In that context one could minimize directly the printing time and at the same time penalize the contact between the desired object and the removable parts, in order to simplify the final detaching operations. %From the numerical point of view, results could be improved by using semi-Lagrangian schemes instead of finite difference schemes, although the CPU time will increase considerably. 

\section*{Acknowledgements}
\REV{Authors want to thank Maurizio Falcone for the useful discussions about the model developed in this paper.}

\appendix

\section{Proof of Lemma \ref{lemma}}
\noindent \emph{Proof.} First of all, it is useful to recall a basic property of the $\tr$ operator. Let us consider two $n\times n$ matrices $\mathbf A$ and $\mathbf B$, and define $\mathbf C:=\mathbf A\mathbf B$. Denote by $\mathbf a_{i,\cdot}$ the $i$-th row of $\mathbf A$ and by $\mathbf b_{\cdot, j}$ the $j$-th column of $\mathbf B$. We have
\begin{equation}\label{trick}
\tr(\mathbf C)\stackrel{\text{def}}{=}\sum_i c_{i,i}=\sum_i \mathbf a_{i,\cdot} \cdot \mathbf b_{\cdot, i}^\intercal=\sum_i\sum_ja_{i,j}\ \! b_{j,i}.
\end{equation}
The Lemma is proved as follows:
\begin{multline*}
|\Grad u|\textup{div}\left(\frac{\Grad u}{|\Grad u|}\right)=
|\Grad u|\sum_{i=1}^n \partial_i\left(\frac{\partial_i u}{|\Grad u|}\right)=\\
|\Grad u|\frac{1}{|\Grad u|^2}\sum_{i=1}^n\left(\partial_i^2 u |\Grad u|-\partial_i u \frac{1}{|\Grad u|}\sum_{j=1}^n\partial_j u\ \partial_i\partial_j u\right)=\\
\triangle u-\frac{1}{|\Grad u|^2}\sum_{i=1}^n \partial_i u\sum_{j=1}^n\partial_j u\ \partial_i\partial_j u=\\
\triangle u-\frac{1}{|\Grad u|^2}\sum_{i=1}^n\sum_{j=1}^n (\Grad u \otimes \Grad u)_{i,j}\ \partial_i\partial_j u=\\
\tr(\GH u)-\sum_{i=1}^n\sum_{j=1}^n \frac{(\Grad u \otimes \Grad u)_{i,j}}{|\Grad u|^2}\ \partial_j\partial_i u
\stackrel{\eqref{trick}}{=}\\
\tr(\GH u)-\tr \left(\frac{\Grad u \otimes \Grad u}{|\Grad u|^2}\GH u\right)=\\
\tr\left(\left(I-\frac{\Grad u \otimes \Grad u}{|\Grad u|^2}\right)\GH u\right).
\qed
\end{multline*}

\section*{References}
%\bibliographystyle{elsarticle-num} 
%\nocite{*}
\bibliography{biblio}

\end{document}